\documentclass{scrartcl}
\usepackage[british]{babel}
\usepackage[textwidth=150mm,top=2cm]{geometry}
\usepackage{dsfont, amsthm, mathtools,  enumerate,amsmath,amsrefs,mathrsfs,lmodern, slantsc}
\usepackage{libertine}

\title{On \textsc{Bernstein}-type theorems}
\author{Peter \textsc{Lewintan}\footnote{peter.lewintan@uni-due.de, University of Duisburg-Essen, Germany.}}
\date{July 18, 2013}

\begin{document}

\maketitle
\begin{center}
\emph{In memory of Franki Dillen}
\end{center}
\thispagestyle{empty}
%\begin{tikzpicture}[remember picture, overlay]
% \node [xshift=1cm,yshift=15cm,rotate=90] at (current page.south west)
% {
% erschienen in: \emph{Pure and Applied Differential Geometry - PADGE 2012},
%In Memory of Franki Dillen, Shaker Verlag 2013,
%168\,--\,174.
% };
%\end{tikzpicture}

\begin{abstract}
We summarize results concerning the \textsc{Bernstein} property of differential equations.
\end{abstract}

In this short overview we will look at entire solutions of partial differential equations of second order. We say a solution to be \emph{entire} if it is defined over the entire plane (${\mathds R}^2$) or over the entire space (${\mathds R}^n$).

\hfill

As we will see, some differential equations possess \emph{only} linear functions as entire solutions, i.e., in these cases the linearity of an entire solution follows from its mere existence, without any boundedness conditions. If a partial differential equation has only affine linear functions as entire solutions, we say that it has the \textsl{\textsc{Bernstein}} \emph{property}, according to a celebrated result of S. N. \textsc{Bernstein} \cite{Lewintan:Bernstein}, which states that every $C^2$-solution of the minimal surface equation
\begin{equation*}
(1+{u_y}^2)u_{xx}-2 u_x u_y u_{xy} + (1+{u_x}^2)u_{yy}=0
\end{equation*}
over the entire plane ${\mathds R}^2$ is necessarily affine linear.

\hfill

We start with the following operator introduced in \cite{Lewintan:Zorina}:
\begin{equation*}
L_{\gamma,\varepsilon}[u]\mathrel{\mathop:}=\left(2\varepsilon + (\gamma+1){u_x}^2+(\gamma-1){u_y}^2\right)u_{xx}+4 u_x u_y u_{xy}+\left(2\varepsilon + (\gamma-1){u_x}^2+ (\gamma +1){u_y}^2\right)u_{yy}
\end{equation*}
with $\gamma, \varepsilon\in{\mathds R}$ and consider the equation
\begin{equation*}
L_{\gamma,\varepsilon}[u]= 0 \qquad \text{over ${\mathds R}^2$}.
\end{equation*}
Without loss of generality we can choose $\varepsilon\in\{-1; 0; 1\}$, for we can obtain entire $C^2$-solutions of $L_{\gamma,\varepsilon}[u]=0$ (with $\varepsilon\neq 0$) via an appropriate scaling of the solutions of $L_{\gamma,\pm 1}[u]=0$ and vice versa.

\hfill

Our equation $L_{\gamma,\varepsilon}[u]=0$ is elliptic if $\varepsilon\gamma>0$ and $\left|\gamma\right|\geq1$. We start with this case and consider other cases later:

\hfill

First we choose $\gamma=\varepsilon=-1$, so that $L_{-1,-1}[u]=0$ corresponds to the familiar minimal surface equation over ${\mathds R}^2$. As we have already mentioned, it has the \textsc{Bernstein} property. The extension of this result to higher dimensions is well-known:

\hfill

The \textsc{Bernstein} theorem was extended to $n\leq 7$, i.e., each entire $C^2$-solution of the minimal surface equation
\begin{equation*}
\operatorname{div}\left(\frac{Du}{\sqrt{1+\left|Du\right|^2}}\right)=0 \qquad \text{over ${\mathds R}^n$}
\end{equation*}
has to be affine linear, cf. \cite{Lewintan:deGiorgi} for $n=3$, \cite{Lewintan:Almgren} for $n=4$ and \cite{Lewintan:Simons} for $n\leq 7$.

\hfill

Surprisingly, the \textsc{Bernstein} theorem fails for dimensions $n\geq 8$ as there exist entire non-linear solutions to the minimal surface equation, cf. \cite{Lewintan:BdGG}. However, under suitable growth conditions on the solution $u$ or its gradient $Du$, one can prove \textsc{Bernstein}-type theorems in every dimension $n\in\mathds N$, cf. e.g. \cite{Lewintan:BdGM}, \cite{Lewintan:Moser}, \cite{Lewintan:CNS}, \cite{Lewintan:Nitsche}, \cite{Lewintan:EH}.

\hfill

Let us now turn to higher codimension $k>1$, so that instead of the minimal surface equation, we consider a system of partial differential equations, the so-called \emph{minimal surface system}. Already in the simplest case of dimension $n=2$ and codimension $k=2$ we obtain several entire non-linear solutions $f\in C^2({\mathds R}^2,{\mathds R}^2)$ of the corresponding system
\begin{equation*}
(1+\left|f_y\right|^2) f_{xx}-2 f_x {\text{\scriptsize $\bullet$}} f_y f_{xy} + (1+\left|f_x\right|^2) f_{yy}=0,
\end{equation*}
for every holomorphic function $f: \mathds C \to \mathds C$, regarded as a map ${\mathds R}^2 \to {\mathds R}^2$, solves it. If $k>1$, the boundedness of the gradient is a sufficient condition for a \textsc{Bernstein}-type theorem only when $n\leq 3$, cf. \cite{Lewintan:CO} and \cite{Lewintan:Fischer}. An analogous result for $n\geq 4$ is wrong, as follows from the \textsc{Lawson-Osserman} cone \cite{Lewintan:LO}\footnote{This cone is an example for a non-analytic \textsc{Lipschitz} solution in higher codimensions.}. Under sufficiently strong assumptions, one can still achieve the linearity of solutions, cf. e.g. \cite{Lewintan:HJW}, \cite{Lewintan:JX}, \cite{Lewintan:Wang}, \cite{Lewintan:JXY}.

\hfill

Returning to the original \textsc{Bernstein} theorem, we can extend it to further differential equations. Such classes of elliptic differential equations (over ${\mathds R}^2$), entire $C^2$-solutions of which are necessarily affine linear, were given e.g. in \cite{Lewintan:Bers}, \cite{Lewintan:Finn}, \cite{Lewintan:Jenkins}, \cite{Lewintan:Simon1}, \cite{Lewintan:Simon2}. The minimal surface equation is included in all these classes. Equations of "minimal surface type" possess the \textsc{Bernstein} property for $n\leq 7$, cf. \cite{Lewintan:Simon1}, and we need additional growth conditions in other dimensions, cf. \cite{Lewintan:Winklmann}.

\hfill

For an elaborated account of the minimal surface case we refer to the monograph \cite{Lewintan:DHT}.

\hfill

Now we choose $\gamma=\varepsilon=1$, so that $L_{1,1}[u]=0$ corresponds to the "wrong minimal surface equation"
\begin{equation*}
(1+{u_x}^2)u_{xx}+2u_x u_y u_{xy}+(1+{u_y}^2)u_{yy}=0.
\end{equation*}
In \cite{Lewintan:Simon2} \textsc{Simon} posed the question whether this equation has the \textsc{Bernstein} property. We can answer in two different ways:
\begin{itemize}
	\item By the separation ansatz $u(x,y)=g(x)+h(y)$ we construct entire non-linear $C^2$-solutions of this equation explicitly. In a similar manner we can determine further (not necessarily elliptic) differential equations without the \textsc{Bernstein} property, cf. \cite{Lewintan:Lewintan}.
	\item We use an explicit criterion of J. C. C. \textsc{Nitsche} and J. A. \textsc{Nitsche} \cite{Lewintan:NN}, which ensures the existence of entire non-linear $C^2$-solutions of certain elliptic differential equations including the wrong minimal surface equation.
\end{itemize}

The \textsc{Nitsche} criterion states:

\hfill

The \textsc{Euler-Lagrange} equation arising from the regular variational integral
\begin{equation*}
\underset{{\mathds R}^2}{\int}F(\left|Du\right|^2)\mathrm d x
\end{equation*}
has entire non-linear $C^2$-solutions if the integral
\begin{equation*}
\overset{\hspace{2mm}\infty}{\underset{\hspace{-3mm}1}{\int}}\frac{1+w\lambda(w)}{2+w\lambda(w)}\cdot\frac{{\ \mathrm d} w}{w}
\qquad \text{with}\qquad \lambda(w)\mathrel{\mathop:}=2\,\frac{F''(w)}{F'(w)}
\end{equation*}
diverges.

\hfill

With this criterion we can treat all the other combinations of $\gamma$ and $\varepsilon$ in the elliptic case, more precisely:

\hfill

In the elliptic case ($\varepsilon\gamma>0$ and $|\gamma|\geq1$) we obtain our equation $L_{\gamma, \varepsilon}[u]=0$ as \textsc{Euler-Lagrange} equation of the functional
\begin{equation*}
\mathscr F_{\gamma,\varepsilon}(u)\mathrel{\mathop:}={\underset{{\mathds R}^2}{\int}}F_{\gamma,\varepsilon}(\left|Du\right|^2)\mathrm d x
\end{equation*}
by setting $w\mathrel{\mathop:}=\left|Du\right|^2$ and
\begin{align*}
F_{\gamma,\varepsilon}(w)\mathrel{\mathop:}=
\begin{cases}\displaystyle
(2|\varepsilon|+|\gamma-1|w)^{\frac{\gamma}{\gamma-1}},&\text{for ($\gamma>1, \varepsilon>0$) or ($\gamma\leq-1, \varepsilon<0$),}\\
\displaystyle\mathrm e^{\frac{w}{2\varepsilon}},&\text{for ($\gamma=1, \varepsilon>0$)}.
\end{cases}
\end{align*}
In all these cases we gain
\begin{equation*}
\lambda_{\gamma,\varepsilon}(w)=2\,\frac{F''_{\gamma,\varepsilon}(w)}{F'_{\gamma,\varepsilon}(w)}=\frac{2}{2\varepsilon + (\gamma - 1)w}.
\end{equation*}
The integral
\begin{equation*}
\overset{\hspace{2mm}\infty}{\underset{\hspace{-3mm}1}{\int}}\frac{1+w\lambda_{\gamma,\varepsilon}(w)}{2+w\lambda_{\gamma,\varepsilon}(w)}\cdot\frac{{\ \mathrm d} w}{w}=\frac{1}{2}\overset{\hspace{2mm}\infty}{\underset{\hspace{-3mm}1}{\int}}\left(\frac{1}{2\varepsilon +\gamma w}+ \frac{1}{w}\right){\ \mathrm d} w
\end{equation*}
diverges for all admissible $\gamma\neq -1$. By the \textsc{Nitsche} criterion, in all these cases there exist entire non-linear $C^2$-solutions of the corresponding \textsc{Euler-Lagrange} equation, i.e. for $\gamma\geq 1, \varepsilon>0$ and $\gamma<-1, \varepsilon<0$ resp. the equation $L_{\gamma, \varepsilon}[u]=0$ does not have the \textsc{Bernstein} property.

\hfill

If we now let $\varepsilon$ tend to zero with $|\gamma|>1$, the integrand converges to
\begin{align*}
F_{\gamma,0}(w)=(|\gamma-1|w)^{\frac{\gamma}{\gamma-1}}.
\end{align*}
By the substitution $\displaystyle p = \frac{2\gamma}{\gamma - 1}$ (for $\gamma\neq 1$) we obtain
\begin{equation*}
F_{\gamma,0}(\left|Du\right|^2)= c(p)\,\frac{1}{p}\,\left|Du\right|^p.
\end{equation*}
Thus, we can associate our functional $\mathscr F_{\gamma,0}$ with the functional
\begin{equation*}
\mathscr F_p(u)\mathrel{\mathop:}=\frac{1}{p}{\underset{{\mathds R}^2}{\int}}\left|Du\right|^p\mathrm d x \qquad \text{where $\displaystyle p=\frac{2\gamma}{\gamma-1}$.}
\end{equation*}
The minimizers of the latter functional are the so-called \emph{$p$-harmonic functions}, the solutions of
\begin{equation*}
\Delta_p u \mathrel{\mathop:}= \operatorname{div}(\left|Du\right|^{p-2}Du) = 0 \quad \text{over ${\mathds R}^2$.}
\end{equation*}

We can regard solutions of $L_{\gamma,0}[u]=0$ as solutions of $\Delta_p u =0$. It is common to introduce the $p$-harmonic functions in the weak sense and not as solutions of $L_{\gamma,0}[u]=0$. As far as the author is aware, it is not clear whether entire $C^2$-solutions of $L_{\gamma,0}[u]=0$ with  $|\gamma|>1$ are necessarily affine linear. However, under suitable growth conditions, cf. \cite{Lewintan:KSZ}, each entire $p$-harmonic function is affine linear. This statement also holds true in higher dimensions.

\hfill

For $\gamma = -1$ we get $p = 1$ and the equation $L_{-1,0}[u]=0$, i.e.
\begin{equation*}
{u_y}^2 u_{xx}- 2 u_x u_y u_{xy} + {u_x}^2 u_{yy}=0,
\end{equation*}
which corresponds to
\begin{equation*}
\Delta_1 u\mathrel{\mathop:}=\operatorname{div} \left(\frac{Du}{\left|Du\right|}\right) = 0.
\end{equation*}
The equation $L_{-1,0}[u]=0$ does not have the \textsc{Bernstein} property, for there are solutions of the form $u(x,y)=g(x)$, with an arbitrary $g\in C^2({\mathds R}, {\mathds R})$, or $u(x,y)={\mathrm e}^{x+y}$.

\hfill

In a more interesting case $\gamma$ tends to $1$ and $p$ tends to $+\infty$. Indeed, the equation $L_{1,0}[u]=0$ corresponds to the equation of the so-called \emph{$\infty$-harmonic function} over ${\mathds R}^2$
\begin{equation*}
{u_x}^2u_{xx} +2 u_x u_y u_{xy} + {u_y}^2u_{yy}=0.
\end{equation*}

The latter has the \textsc{Bernstein} property, cf. \cite{Lewintan:Aronsson}. One can extend this result to higher dimensions if additional regularity be assumed,  more precisely:

\hfill

Each entire $C^4$-solution of
\begin{equation*}
\Delta_\infty u \mathrel{\mathop:}= \underset{j,k=1}{\overset{n}{\sum}}u_{x_j} u_{x_k} u_{x_j x_k}=0 \qquad\text{over ${\mathds R}^n$}
\end{equation*}
is necessarily affine linear, cf. \cite{Lewintan:Yu}. It is not clear whether $C^2$-regularity suffices here. If so, the equation $\Delta_\infty u=0$ would have the \textsc{Bernstein} property in all dimensions, in contrast to the minimal surface equation.

\hfill

Concerning all the other combinations of $\varepsilon$ and $\gamma$, we state the following:

\hfill

For $\varepsilon=\gamma = 0$ we have entire non-linear $C^2$-solutions:\quad $u(x,y)=x^2+y^2$ solves  $L_{0,0}[u]=0$.

\hfill

Also, for $\varepsilon > 0$ and $\gamma = -1$ our equation does not have the \textsc{Bernstein} property, for $L_{-1,1}[u]=0$ admits solutions of the form $u(x,y)= x + h(y)$, with an arbitrary $h\in C^2(\mathds R,\mathds R)$. However, $L_{-1,1}[u]=0$ corresponds to the maximal surface equation
\begin{equation*}
(1-{u_y}^2)u_{xx}+2u_xu_yu_{xy}+(1-{u_x}^2)u_{yy}=0,
\end{equation*}
and under constraint $|Du|_{C^0}<1$ its entire solutions are affine linear. This is also valid in higher dimensions, cf. \cite{Lewintan:Calabi} and \cite{Lewintan:CY}.

\hfill

For $\varepsilon < 0$ and $\gamma = 1$ our equation $L_{\gamma,\varepsilon}[u] = 0$ arises in the study of isentropic irrotational steady plane flows, cf. \cite{Lewintan:CF}.

\end{document}